\documentclass[10pt]{amsart}
\usepackage{amsmath,amssymb,graphicx}

\title{Multiplication maps of linear systems on smooth projective toric surfaces}
\author{Najmuddin Fakhruddin}

\address{School of Mathematics, Tata Institute of Fundamental Research, Homi Bhabha Road, Mumbai 400005, INDIA}
\email{naf@math.tifr.res.in}
\renewcommand{\L}{\mathcal{L}}
\newcommand{\M}{\mathcal{M}}
\renewcommand{\O}{\mathcal{O}}
\newcommand{\R}{\mathbf{R}}
\newcommand{\Z}{\mathbf{Z}}
\newcommand{\la}{\langle}
\newcommand{\ra}{\rangle}
\newcommand{\s}{\sigma}
\renewcommand{\S}{\Sigma}
\newcommand{\bs}{\backslash}

\newtheorem{thm}{Theorem}
\newtheorem{lem}{Lemma}

\theoremstyle{remark}

\begin{document}

\begin{abstract}
Let $X$ be a smooth projective toric surface and $\L$ and $\M$ two line bundles
on $X$. If $\L$ is ample and $\M$ is generated by global sections, then we show
that the natural map $H^0(X,\L) \otimes H^0(X,\M) \to H^0(X, \L \otimes \M)$ is 
surjective. 
We also consider a generalization to the case when $\M$ is
arbitrary line bundle with $ h^0(X,\M) > 0$.
\end{abstract}

\maketitle

In this note we shall prove the following results:

\begin{thm}
  Let $X$ be a smooth projective toric surface, $\L$ an ample line bundle on
  $X$, and $\M$ a line bundle on $X$ which is generated by global
  sections. Then the multiplication map $H^0(X,\L) \otimes H^0(X, \M)
  \to H^0(X, \L \otimes \M)$ is surjective. \label{thm:sections}
\end{thm}

\begin{thm}
  Let $X$ be a smooth projective toric surface and $\L$ an ample line bundle
  on $X$. Then there exists a constant $C(\L)$ such that for all line
  bundles $\M$ on $X$ with $H^0(X, \M) \neq 0$,
  $dim(coker[H^0(X,\L) \otimes H^0(X, \M) \to H^0(X, \L \otimes \M)])
  \leq C(\L)$. \label{thm:effective}
\end{thm}

We do not know if either of these results holds if $dim(X) > 2$ or if $dim(X)=2$
and $X$ is singular. Similar questions have been raised by Oda \cite{oda-problems}.

\section{}
We refer the reader to  \cite{fulton-toric}
or \cite{oda-toric} for basic facts about toric varieties.

From now on $X= X(\Delta)$ will always be a smooth projective toric
surface associated to a fan $\Delta$ in $N_{\R}$, where $M$ (resp.
$N$) is the lattice of characters (resp. co-characters) of a
$2$-dimensional algebraic torus. Any divisor on $X$ is linearly
equivalent to a divisor $D= \sum_{i=1}^n a_i D_i$ with $a_i \in \Z$,
and the $D_i$'s the divisors on $X$ invariant under the torus action;
these are in $1$-$1$ correspondence with the rays in $\Delta$.  (We
assume that the indices are chosen so that $D_i$ and $D_{i+1}$
correspond to rays forming the boundary of a cone in $\Delta$.  Here,
and in what follows, we assume that subscripts are considered modulo
$n$.).  Moreover, if $D$ is effective then we may assume that all $a_i
\geq 0$ and we denote by $P_D$ the corresponding polygon in $M_{\R}$.
We let $v_i$ be the minimal lattice vector in the ray corresponding to
$D_i$.

The proof of Theorem \ref{thm:sections} is easily reduced to
combinatorial statements about convex polygons, using the well-known
dictionary relating equivariant divisors $D$ with $\O(D)$ generated by
sections and lattice polygons. We may assume that $\L$ and $\M$ are of
the form $\O(D)$ and $\O(E)$ for some equivariant divisors $D=
\sum_{i=1}^n a_i D_i$ and ${E} = \sum_{i=1}^n b_i D_i$.  Then it
suffices to prove that $(P_D \cap M) + (P_{E} \cap M) = (P_D + P_{E})
\cap M$.

For a divisor $D$ as above with $\O(D)$ generated by sections,
let $\s_i(D) = P_D \cap L_i$ where
 $L_i = \{ u \in M_{\R} |  \la u,v_i\ra = - a_i \}$. If $D$
is ample then $\s_i(D)$ is always an edge of $P_D$ but in general it could
also be a vertex.

\begin{lem}
Let $D$, $E$ be as above and assume that $D$ is ample and $\O(E)$
is generated by sections. Then $\s_i(D + E) = \s_i(D) + \s_i(E)$
for all $i \in [1,n]$.  \label{lem:edges}
\end{lem}
\begin{proof}
This follows easily from the fact that $P_{D+E} = P_D + P_E$.
\end{proof}

The proof of the following lemma is also left to the reader.
\begin{lem}
Let $[a_1,b_1]$ and $[a_2,b_2]$ be closed intervals of $\R$.
Suppose $[a_1,b_1] \cap \Z \neq \emptyset $ and $a_2,b_2 \in \Z$.
Then any element $z$ of $[a_1 + a_2, b_1 + b_2] \cap \Z$ is
of the form $c_1 + c_2$ for some $c_i  \in [a_i,b_i] \cap \Z$, $i=1,2$.
\label{lem:trivial}
\end{lem}

To prove Theorem \ref{thm:sections} we will first reduce to the case where
$P_{E}$ is a triangle of a special kind, and then explicitly
 prove the equality in this case.

\begin{proof}[Proof of Theorem \ref{thm:sections}]
  We first dispose off the trivial cases: If $P_{E}$ is a point then
  the statement is obvious and if $P_{E}$ is $1$-dimensional, hence a
  line segment, then the proof is elementary and we leave it to the
  reader (use Lemma \ref{lem:trivial}).
  
  Let $s: M_{\R} \times M_{\R} \to M_{\R}$ be the sum map $(x,y)
  \mapsto x+ y$.  Then $P_{D + {E}} = P_D + P_{E} = s(P_D \times
  P_{E})$. Let $p \in P_{D + {E}} \cap M$ and let $Q = s^{-1}(p) \cap
  (P_D \times P_{E})$. This is a convex polygon  (possibly degenerate) 
  in $ M_{\R} \times
  M_{\R}$ and we let $Q_i = \pi_i(Q)$, where $\pi_i$, $i= 1,2$, are
  the two projections. So $Q_1 \subset P_D$ and $Q_2 \subset P_{E}$
  are also convex polygons.
  
  Let $(q_1,q_2) \in Q$ be such that $q_2$ is a vertex of $Q_2$. If
  $q_2$ is in the interior of $P_E$ then $q_1 \in Q_1$ must be a
  vertex of $P_D$ (sic), hence $q_1 \in M$. Then $q_2 \in M$ and $p =
  q_1+ q_2$, so we are done.  Otherwise, since $Q_2$ must have at
  least one vertex, it follows that there exists a point $q \in Q_2$ which
  lies on the boundary of $P_{E}$. If $q \in M$, then we are done so
  we may assume that $q$ lies in the interior of an edge $\s$ of
  $P_{E}$. We let $m_1, m_2$ be the two end points of $\s$.
  
  Recall that $P_{E} = \{ u \in M_{\R} | \la u,v_i\ra \geq - b_i
  \mbox{ for all } i \}$.  We may assume that the edge $\s$
  corresponds to $v_1$, so $\la m_i, v_1\ra = -b_1$, $i=1,2$.  Now for
  each $i = 1, \ldots, n$, let $c_i = min\{ c \in \Z | \la m_j, v_1\ra
  \geq -c \mbox{ for } j = 1,2\}$.  Then $c_1 = b_1$ and $c_i \leq
  b_i$ for all $i$.  Let $P = \{ u \in M_{\R} | \la u,v_i\ra \geq -
  c_i \mbox{ for all } i \}$; by construction $P \subset P_{E}$ and
  $m_1$ and $m_2$ are vertices of $P$.  If $P = \s$ then we are
  done, so we may assume that $P$ is $2$-dimensional.  Without loss of
  generality we may assume that $\la m_1, v_2\ra > \la m_2,v_2\ra $
  and we let $k = max \{ i \in [2,n] | \la m_1, v_i\ra > \la
  m_2,v_i\ra \}$. So for $i \in [1,k]$, $c_i = - \la m_2,v_i\ra $ and
  for $i \in [k+1, n]$, $c_i = - \la m_1,v_i\ra $.  By our assumption
  that $P$ is not $1$-dimensional it follows that $\la m_1, v_i\ra
  \neq \la m_2,v_i\ra $ for all $i \in [2,n]$. Then $P = \{ u \in
  M_{\R} | \la u,v_i\ra \geq - c_i \mbox{ for } i=1,k, k+1 \}$, hence
  $P$ a triangle. Since $X$ is smooth it follows that the third vertex
  is also in $M$, so $P$ corresponds to an equivariant divisor on $X$
  whose associated line bundle is generated by sections.
  
  Since $q \in P$ by construction, by replacing $P_{E}$ by $P$ we have
  reduced to the case when $P_E$ is a triangle with the further
  property that there exists an $i \in [1,n]$ such that $\s_i(E)$ and
  $\s_{i+1}(E)$ are both (non-degenerate) edges of $P_E$.  By using the
  basis of $M$ dual to $\{v_i,v_{i+1}\}$, and after possible translation
  by elements of $M$ (which does not affect the hypotheses or the
  conclusion), we have the following picture: $P_E$ is the convex span
  of the points $(0,0), (a,0), (0,b)$, for some $a,b > 0$, $P_D$ is
  entirely contained in the first quadrant, and $(0,0)$ is also a
  vertex of $P_D$ (consequently $P_D$ must also have edges along the
  positive $x$ and $y$ axes).
  
  We shall now complete the proof of the theorem by analysing this
  case.  Decompose the region $P_D + P_E \bs P_D$ as a union
  of the three regions, $A, B$ and $C$, as illustrated in the 
  figure below
  --- to see that this is correct we use Lemma \ref{lem:edges}.
  Note that $A$ or $B$ may be empty; this happens precisely when
  $k=2$ or $k=n-1$.

  We claim that any lattice point in
  the region $A$ is of the form $m + (x,0)$ where $m \in P_D \cap M$
  and $0 \leq x \leq a$.  This is because the trapezium, two of whose
  sides are the base of the triangle $P$ and the edge $U$ of $P_D$, is
  contained in $P_D$ and both these sides contain at least two lattice
  points each. Thus  each horizontal line which contains a lattice element
  of $A$ also contains a lattice element of $P_D$, so the claim
  follows by Lemma \ref{lem:trivial}.
  By a symmetric argument,
  any lattice point in the region $B$ is of the form 
  $m + (0,y)$ where $m \in P_D \cap M$ and $0 \leq y \leq b$.

  Any point in the region $C$ is
  contained in $P +P_E$. Since $P$ and $P_E$ are similar triangles
  (i.e. are translates of multiples of the same triangle) one easily
  sees that any lattice point in $P + P_E$ is the sum of lattice
  points in $P$ and $P_E$.

\vspace{1cm}
\begin{center}
\includegraphics[width=5in]{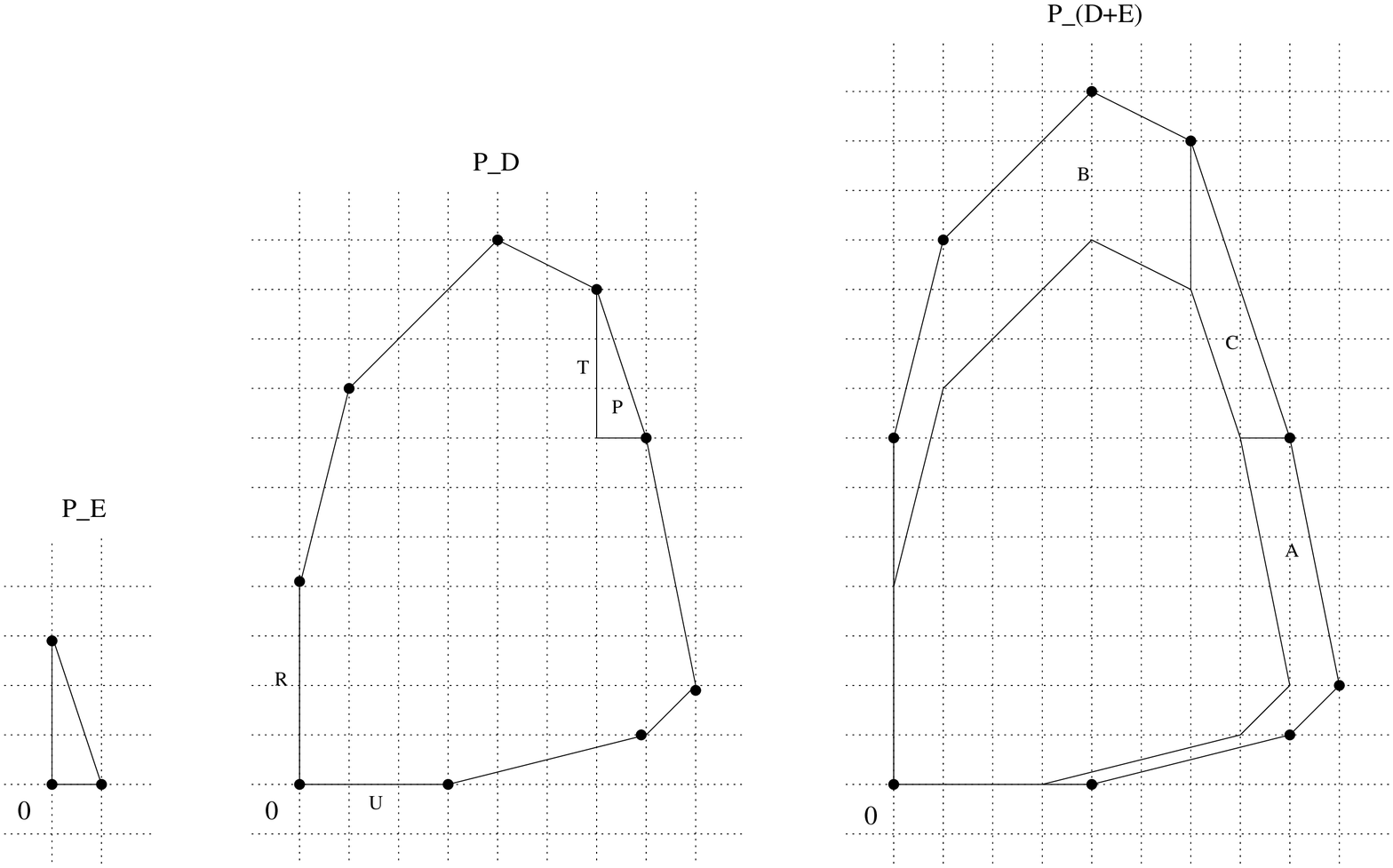}
\end{center}
\end{proof}

\section{}

The following lemma is the key to the deduction of Theorem
\ref{thm:effective} from Theorem \ref{thm:sections}.

\begin{lem}
  Let $D = \sum_{i=1}^n a_i D_i$ be an effective divisor on $X$. Then
  there exists integers $b_i$, $0 \leq b_i \leq a_i$, such that for
  $D' = \sum_{i=1}^n b_i D_i$, $\O(D')$ is generated by its sections
  and the natural map $H^0(X, \O(D')) \to H^0(X, \O(D))$ is an
  isomorphism. \label{lem:reduction}
\end{lem}

\begin{proof}
  Let $P_D$ be the polygon associated to $D$. Let $S = P_D \cap M$ and
  let $P$ be the convex hull of the points in $S$.  For each $i$, $1
  \leq i \leq n$, let $b_i = min\{c \in \Z | \la s,v_i\ra \geq -c
  \mbox{ for all } s \in S \}$.  Since $0 \in S$ and $ S \subset P_D$,
  it follows that $0 \leq b_i \leq a_i$.  We claim that $P$ is equal
  to $P' := \{ u \in M_{\R} | \la u,v_i\ra \geq - b_i \mbox{ for all }
  i \}$.
  
  Consider the lines $L_i = \{ u \in M_{\R} | \la u,v_i\ra = - b_i
  \}$.  By construction, each of these lines contains a point $s_i$ of
  $S$ and all points of $S$ are contained in the ``positive'' half
  planes $H_i =\{ u \in M_{\R} | \la u,v_i\ra \geq - b_i \}$.  It
  suffices to show that the point $m_i = L_i \cap L_{i+1}$ is in $S$
  for all $i$. Suppose not; then there exists $i$ and $j$ such that
  $\la m_i, v_j\ra < -b_j$. Since $\la s_i, v_j\ra \geq -b_j$ and $\la
  s_{i+1}, v_j\ra \geq -b_j$, it follows $L_j$ must intersect $L_i$
  and $L_{i+1}$ along the segments $[s_i, m_i]$ and $[s_{i+1},m_i]$.
  But looking at the normal rays, this says that the ray along $v_j$
  must lie between the rays along $v_i$ and $v_{i+1}$. But we have
  assumed that the $v_i$'s are cyclically ordered, so this is a
  contradiction.
  
  Since $P=P'$, it follows that $P$ corresponds to a divisor $D'$ on
  $X$ of the required form.
\end{proof}

\begin{lem}
  There exists a constant $C$ depending only on $X$ with the following
  property: Let $D, D'$ be as in Lemma \ref{lem:reduction} and let $J
  = \{ i \in [1,n] | b_i < a_i\}$. Then the number of lattice points
  on all the edges of $P_D'$ whose normal ray contains $v_j$ for some
  $j \in J$ is bounded by $C$. \label{lem:length}
\end{lem}

\begin{proof}
  This follows because there are only finitely many rays in $\Delta$,
  the fan of $X$.  If we choose a Euclidean metric on $M_{\R}$, then
  there are only finitely many possibilities for the angles at the
  vertices of any polygon associated to a divisor on $X$. Bounding the
  number of lattice points is also the same as bounding the lengths of
  edges.  Suppose that the length of the edge corresponding to $v_j$
  is not bounded.  Let $D'' = \sum_{i=1}^n c_i D_i$ with $c_i = b _i$
  for $ i \neq j$ and $c_j = b_j + 1$. Then $P_{D'} \subset P_{D''}
  \subset P_D$ and if the length of is sufficiently large, $P_{D''}$
  would contain lattice points not contained in $P_D$, contradicting
  the defining property of $D'$.
\end{proof}

The next lemma implies that even though there may
be infinitely many divisors $D$ giving rise to the same $D'$ as in
Lemma \ref{lem:reduction}, upto translation by elements of $M$ there
are only finitely many polygons occuring as connected components of
$P_D \bs P_{D'}$, where $D$ ranges over all effective divisors
on $X$.

\begin{lem}
  There exists
  a constant $C_2$ depending only on $X$ such that 
  for $D$, $D'$ and $J$ as in Lemma \ref{lem:length}
  and 
  $J' =  [j_1,j_2]$ any subinterval of $[1,n]$ contained in $J$,
  exactly one of the following holds: \\
  1) $[1,n] \bs J'$ contains at most one element. In this case
  all the edges of $P_{D'}$ have length $\leq C_1$, and there exists
  $J'' \subset [1,n]$ such that $a_{j} \leq C_2$ for $j \in J''$ and
  such that the polyhedron $P(J'') = \{u \in M_{\R} | \la u, v_{j}\ra
  \geq -a_{j} \mbox{ for all } j \in J''\}$
  is bounded. \\
  2) $[1,n] \bs J'$ contains at least two elements, so $j_1 -1$
  and $j_2 +1$ are distinct elements of $[1,n]$.  
  Consider the lines $L_{j_1 -1}$ and $L_{j_2 + 1}$. Then either \\
  2a) The lines intersect in a point $p$ such that any line segment
  joining $p$ and $\s_j$ for any $j \in J'$ does not contain
  any point of $P_{D'}$ except for an endpoint. Or \\
  2b) The two lines are parallel. Then there exists $j \in J'$
  such that $a_j - b_j \leq C_2$. Or \\
  2c) There exists a
  subset $J''$ of $J$ such that $a_{j} - b_{j} \leq
  C_2$ for $j \in J''$ and such that the
  region
\[
P(J'') = \{u \in M_{\R} | \la u, v_{j}\ra \geq -a_{j} \mbox{ for all }
j \in J'' \cup \{j_1 -1, j_2 +1\}\} 
\]
 is bounded.
\end{lem}

\begin{proof}
  The lemma is esentially a consequence of Lemma \ref{lem:length}. Since
  there are only finitely many possibilities for the lengths of
  the edges of $P_{D'}$ corresponding to $j \in J$, the number 
  of possible configurations
  (upto translation) of the subset of the
  the boundary of $P_{D'}$ which is the union of the edges
  corresponding to the $j$'s in $J'$ is also finite. (In case
  1), even the number of possible $D'$ is finite.) So it is enough
  to find a constant $C$ which works in each case separately,
  since we can then let $C_2$ be the maximum of all these.

  First assume that we are in case 1). Let $\S(J')$ be the
  collection of subsets $J''$ of $J'$ such that the rays 
  of $\Delta(X)$ corresponding to $j \in J''$ give rise
  to a complete fan i.e. any open half-space in $N_{\R}$
  must contain one of these rays; so these are precisely
  the susbsets $J''$ for which $P(J'')$ is always bounded.
  Suppose the conclusions in case 1) do not hold.
  Since there are only a finite number of possible $D'$,
  we may consider each of them separately, so we may assume
  that there is no constant which works for some fixed $D'$.
  Since $J'$ is a finite set it follows that there exists
  a sequence of divisors $D^l = \sum_{i=1}^n a_i^lD_i$
  with ${D^l}' = D'$ and a subset $J'''$ of $J'$ such that
  $a_i^l \to \infty$ as $j \to \infty$ for
  all $i \in J'''$. Furthermore $J''' \cap J'' \neq \emptyset$
  for all $J'' \in \S(J')$. It follows that 
  $\cup_l P_{D^l} \supset \{u \in M_{\R} | \la u, v_i \ra \geq -b_i \mbox{ for all } i \in J' \bs J'''\}$.
  But $J' \bs J'''$ is not in $\S(J')$ so $\cup_l P_{D^l}$ contains
  an unbounded polyhedron and hence must contain infinitely many elements
  of $M$. But this contradicts the assumption that ${D^l}' = D'$ for all $l$.


  Case 2) is handled in an analogous manner, the remarks at the beginning
  of the proof allowing us to consider essentially one $D'$ at a time.
  For 2a) there is nothing to prove and 2b) is elementary.
  For 2c) we let $S(J')$ be as above except that we require
  that $J'' \cup \{j_1,j_2\}$ give rise to a complete fan; it follows
   by assumption that  $S(J') \neq \emptyset$.
  The reason for the $a_j - b_j$ occuring here, instead of just the
  $a_j$ in case 1), is because we only have finiteness of possible
  configurations upto translation. (Note that 
  the lengths of the edges correponding to $i \notin J$
  have no effect, since the claim is ``local'' around a given $J'$.)
\end{proof}

\begin{proof}[Proof of Theorem \ref{thm:effective}]
  Let $\L = \O(D)$ and $\M = \O(E)$. Let $E'$ be the divisor
  associated to $E$ using Lemma \ref{lem:length}. By Theorem
  \ref{thm:sections} it follows that the map $H^0(X, \O(D)) \otimes
  H^0(X, \O(E')) \to H^0(X, O(D) \otimes \O(E'))$ is surjective or
  equivalently the map $(P_D \cap M) + (P_{E'} \cap M) \to (P_{D + E'}
  \cap M)$ is surjective. By combining the previous three lemmas, it
  follows that that there are only finitely many possibilities for the
  connected components of $P_{D+E} \bs P_{D + E'}$, upto
  translation by lattice points. This is because $D$ is fixed, so the
  lengths of the edges of $P_{D+E'}$ corresponding to $j \in J$ are
  bounded independently of $E$ (use Lemma \ref{lem:edges}).  The
  number of lattice points in $P_{D+E} \bs P_{D + E'}$ can thus
  be bounded by a constant depending only on $D$, whence the theorem.
\end{proof}


\emph{Acknowledgements.} This note was written in response to a question of V.~Srinivas.
I thank W.~Fulton for informing me about \cite{oda-problems} and T.~Oda for sending it
to me.


\end{document}